\documentclass[12pt]{article}
\usepackage{amsfonts, amsmath, amsthm}
\usepackage[all]{xy}

\newtheorem{theorem}{Theorem}[section]
\newtheorem{lemma}[theorem]{Lemma}
\newtheorem{conj}[theorem]{Conjecture}
\newtheorem{cor}[theorem]{Corollary}
\newtheorem{prop}[theorem]{Proposition}

\def\NN{\mathbb{N}}

\def\QQ{\mathbb{Q}}
\def\RR{\mathbb{R}}
\def\ZZ{\mathbb{Z}}
\def\calO{\mathcal{O}}

\def\alg{\mathrm{alg}}
\def\an{\mathrm{an}}
\def\con{\mathrm{con}}

\def\imm{\mathrm{imm}}
\def\sep{\mathrm{sep}}
\def\perf{\mathrm{perf}}

\def\beq{\begin{equation}}
\def\eeq{\end{equation}}
\def\map#1#2#3{#1\!:\!#2 \to #3}

\def\fp{\frac{1}{p}}

\def\GK{\Gamma^K}
\def\GL{\Gamma^L}
\def\Gperf{\Gamma^{\perf}}
\def\Gsep{\Gamma^{\sep}}
\def\Gimm{\Gamma^{\imm}}
\def\Galg{\Gamma^{\alg}}

\def\Galgcon{\Galg_{\con}}
\def\Gan{\Gamma_{\an}}
\def\Gancon{\Gamma_{\an,\con}}
\def\Gcon{\Gamma_{\con}}
\def\Gimmcon{\Gimm_{\con}}

\def\Gperfcon{\Gperf_{\con}}
\def\Gsepcon{\Gsep_{\con}}

\def\Oan{\Omega_{\an}}

\def\be{\mathbf{e}}
\def\bof{\mathbf{f}}
\def\bv{\mathbf{v}}
\def\bw{\mathbf{w}}

\DeclareMathOperator{\Gal}{Gal}
\DeclareMathOperator{\Hom}{Hom}

\DeclareMathOperator{\Spec}{Spec}

\begin{document}

\title{Unipotency and semistability of overconvergent $F$-crystals}
\author{Kiran S. Kedlaya \\ University of California, Berkeley}
\date{\today}

\maketitle

\begin{abstract}
We introduce the notions of semistability and potential semistability
of overconvergent $F$-crystals over an equal characteristic local field.
We establish their equivalence with the
notions of unipotency and quasi-unipotency given by Crew, and recast 
the conjecture that every overconvergent crystal is quasi-unipotent
in terms of potential semistability.
\end{abstract}

\section{Introduction}

Crystals are the basic objects in any $p$-adic cohomology theory for
schemes of characteristic $p>0$. In particular,
overconvergent $F$-crystals are the principal objects in
Berthelot's theory of rigid cohomology \cite{bib:ber4},
which subsumes earlier constructions of Dwork and Monsky-Washnitzer for
affine schemes, and crystalline cohomology for proper schemes.
(A fuller development of rigid cohomology appears in \cite{bib:ber2}.)
They are also related to Galois representations over a mixed
characteristic discrete valuation ring, by work
initiated by Fontaine and pursued by Fontaine, Colmez and others.

Unfortunately, the global theory of crystals is marred by some gaps
in knowledge in
the local theory (i.e., the theory of crystals over $\Spec k[[t]]$
or $\Spec k((t))$). For example, Crew \cite{bib:crew2} establishes good
structural properties of an overconvergent $F$-crystal on a curve (finiteness
of cohomology plus analogues of some results of Weil II) only under the
local hypothesis that everywhere the crystal is ``quasi-unipotent''.
Crew also suggests that perhaps every
overconvergent crystal over $\Spec k((t))$, or at least every such
crystal ``of geometric origin'',
is quasi-unipotent.

The purpose of this paper is to introduce the notions of semistability
and potential semistability, to establish their equivalence with the
notions of unipotency and quasi-unipotency given by Crew, and to recast 
the conjecture that every overconvergent crystal is quasi-unipotent
in terms of potential semistability. One reason to do
this is that unlike the definition of unipotency, the definition of
semistability avoids referring to any rings which are not discrete
valuation rings. This makes it easier to handle by arguments involving
reduction modulo powers of $p$.

Some of the consequences of this equivalence will be realized in
subsequent papers. These include the facts that de~Jong's extension
theorem \cite{bib:dej1} holds for quasi-unipotent overconvergent $F$-crystals,
and that crystals ``of geometric origin''
are potentially semistable. (The finite dimensionality of rigid cohomology
of varieties follows from this assertion plus results of Crew; it has 
also been established directly
by Berthelot \cite{bib:ber1}.) One mixed characteristic consequence
upon which we will not
dwell further is Berger's proof \cite{berger} that every crystalline
representation
of the Galois group of a local field is of finite height, a result conjectured
by Fontaine.

\section{Some auxiliary rings}

In the next few sections,
 we give detailed constructions of the various coefficient
rings that occur in the local study of crystals. These rings are related to
each other by various augmentations and restrictions; to keep track of these,
we introduce several base rings, which are notated by individual symbols,
and notate the others by attaching
``decorations'' to the names of the base rings. Unfortunately, certain sets of
decorations interact in unexpected ways. We have attempted to flag any such
interactions that may cause trouble, but we make no guarantees about the use of
combinations of decorations not explicitly mentioned in the text.

The development in this section and the next
is largely
modeled on \cite[Section~4]{bib:dej1}. However, some notations has been
changed, some definitions have been made in slightly greater generality,
and some constructions are discussed here in greater detail
for future reference.

\subsection{Liftings to characteristic zero}



First and foremost, fix an algebraically closed field $k$ of characteristic $p$,
and let $W$ be its ring of Witt vectors.
Denote by $\sigma$ the Frobenius map $x \mapsto
x^p$ on $k$ and its canonical lift to $W$. Now let $\calO$
be a finite totally ramified extension of $W$ admitting an automorphism lifting
$\sigma$; we fix a choice of said automorphism and call it $\sigma$ as well.
Let $\pi$ denote a uniformizer of $\calO$ and
$|\cdot|$ the $p$-adic absolute value on $\calO$, normalized so that
$|p| = p^{-1}$. Let $\calO_0$ denote the subset of $\calO$ fixed by $\sigma$.
Because $k$ is algebraically closed, the equation $\lambda^\sigma
= \mu \lambda$ has an equation for any $\mu \in \calO$ with $|\mu| = 1$;
thus every element of $\calO$ can be written as an element of $\calO_0$
times a unit of $\calO$.

Our first main task is to construct complete discrete valuation rings
of characteristic $(0,p)$ lifting the fields in the following tower,
subject to several restrictions. We want these rings to contain $\calO$,
and we want them to admit compatible actions of $\sigma$.
Moreover, if $L = k((u))$ is a finite extension of $K=k((t))$,
we want the tower with the lift of $K$ on the bottom to be compatible with
the tower with the lift of $L$ on the bottom.
\[
\xymatrix{
& K^{\imm} \ar@{-}[d] & \\
& K^{\alg} \ar@{-}[dl] \ar@{-}[dr] & \\
K^{\sep} \ar@{-}[dr] & & K^{\perf} \ar@{-}[dl] \\
& K = k((t)) & \\
}
\]
Here $K^\perf, K^\sep, K^\alg$ denote the perfect,
separable and algebraic closures, respectively, of $K$, while $K^\imm$ denotes
the maximal immediate extension of $K$ in the sense of Kaplansky \cite{bib:kap}. Explicitly,
$k((t))^{\imm}$ consists of generalized power series $\sum_{i \in \QQ} c_i t^i$,
with $c_i \in k$,
for which $c_i = 0$ outside of a well-ordered subset of $\QQ$.

Our first problem is to lift $K$; we cannot use its ring of Witt vectors for
our purposes, because $K$ is not perfect. Instead, we construct its lift $\GK$
as the ring of power series $\sum_{i \in \ZZ} c_i t^i$, with $c_i \in \calO$,
with the property that
for each $\epsilon < 1$, the set of $i \in \ZZ$ such that $|c_i| \geq \epsilon$
is bounded below. (We will drop the decoration $K$ in circumstances where this is
unambiguous.)
This construction singles out a distinct element $t$ of $\GK$
which lifts a uniformizer of $K$, but this distinction is illusory; if $t_1$ is another
residual uniformizer, then each element of $\GK$ has a unique expression as
$\sum_{i \in \ZZ} d_i t_1^i$ such that $\{i \in \ZZ: |d_i| \geq \epsilon\}$ is bounded
below for each $\epsilon < 1$. We will occasionally use the notation $[t^i]x$
to refer to the coefficient of $t^i$ in the expansion of $x$ as a power series in $t$.

Now $\GK$ is clearly what we wanted, a complete discrete valuation ring with
residue field $K$. (Otherwise put, $\GK$ is a Cohen ring of $K$ tensored
over $W$ with $\calO$.)
Let $|\cdot|$ denote the corresponding absolute value, normalized
to be compatible with $|\cdot|$ on $\calO$.
We define a Frobenius $\sigma$ on $\GK$ to be any ring endomorphism of $\GK$
such that $|x^\sigma - x^p| < 1$ for all $x \in \GK$.

\begin{prop} \label{prop:extgam}
Let $R$ be a complete discrete valuation ring of mixed characteristic $(0,p)$,
with residue field $K$, equipped with a Frobenius $\sigma$ lifting the $p$-th power
map on $K$. Let $L$ be an extension of $K$ which is a
separable extension of $K^{1/p^n}$
for some $n \in \NN$.
Then there exists a complete discrete valuation ring $S$ over $R$ with residue field
$L$, such that $\sigma$ extends to an endomorphism of $S$ lifting the $p$-th power
map on $L$. Moreover, this extension is canonical, in that if $R_1, K_1, L_1$ is
an analogous set of rings and there exist maps $R \to R_1$ and $L \to L_1$ compatible
with Frobenius and yielding the same map $K \to K_1$, then there is a unique map
$S \to S_1$ compatible with Frobenius and making the following diagram commute.
\[
\xymatrix{
S \ar@{-->}[rr] \ar@{-}[dd] \ar@{->>}[dr]
& & S_1 \ar@{-}[dd] \ar@{->>}[dr] & \\
& L \ar[rr] \ar@{-}[dd] & & L_1 \ar@{-}[dd] \\
R \ar@{->>}[dr] \ar[rr] & & R_1 \ar@{->>}[dr] & \\
& K \ar[rr] & & K_1
}
\]
\end{prop}
Beware that there is no canonical lift if $L$ is purely inseparable but not of
the form $K^{1/p^n}$, e.g., if $K$ is the fraction field of $k[x,y]$
and $L = K(x^{1/p} + y^{1/p})$.

\begin{proof}
It suffices to consider two cases: $L$ is separable over $K$, or $L = K^{1/p}$.

In case $L$ is normal and separable over $K$, let $P(x)$ be a separable polynomial such that
$L = K[x]/(P(x))$ (which exists by the primitive element theorem).
Choose a lift $\tilde{P}(x)$ of $P$ to $R$, and define $S$
to be $R[x]/(\tilde{P}(x))$. Then $S$ is a complete discrete valuation ring
with residue field $L$, so it is henselian. In particular, if $Q(x)$ is any other
separable polynomial
such that $L = K[x]/(Q(x))$ and $\tilde{Q}(x)$ is any lift of $Q$ to $R$, then 
$\tilde{Q}(x)$ has a root in $S$, so $S$ is well-defined up to isomorphism.

To extend $\sigma$ to $S$, set $x^\sigma$ to be the unique root of the polynomial
$\tilde{P}^\sigma(x)$ congruent to $x^p$ modulo $\pi$ (which again exists by Hensel's
Lemma). It is easily verified that this definition is also independent of the choice
of $P$ and its lift $\tilde{P}$, which is precisely to say that the desired
canonicality holds.

In case $L = K^{1/p}$, take $S$ to be a ring isomorphic to $R$ and let $\phi: R \to S$
denote the isomorphism. Now map $R$ into $S$ by sending $r \in R$ to $\phi(r)^\sigma$.
Canonicality is immediate in this case: if $S_1$ is analogously defined to be isomorphic
to $R_1$ via $\phi_1: R_1 \to S_1$, we map $S$ to $S_1$ by sending $s$ to
$\phi_1(\phi^{-1}(s)^\sigma)$.
\end{proof}

We must still establish that if $K = k((t))$ and $R = \GK$,
then $S$ is isomorphic to $\GL$. This is obvious from the above construction
when $L/K$ is purely inseparable. On the other hand, suppose $L/K$ is separable,
with $L = K(u)$ where $u$ is a root of the polynomial $P(x)$ over $K$. We can construct
an embedding of $\GK$ into $\GL$ by mapping some chosen lift of $t$ to $\GK$
to some other chosen lift of $t$, viewed inside $L$, to $\GL$; for simplicitly
of language, we will regard $\GK$ as a subring of $\GL$ using this embedding.
Then $\GL$ is finite over $\GK$: a finite set of generators of $L/K$
lifts to a set of generators of $\GL$ over $\GK$ by Nakayama's Lemma.
In particular, $\GL$ contains $S$, since it is henselian and so contains a root
of any separable polynomial which has a residual root. Since $S$ and $\GL$ have
the same residual degree, namely $[L:K]$, they coincide.

We lift $K^{\sep}$, $K^{\perf}$, and $K^{\alg}$ by taking $p$-adic completions of
direct limits:
\begin{align*}
\Gsep &= \left( \varinjlim_{\text{$L/K$ separable}} \quad \GL \right)^\wedge \\
\Gperf &= \left( \varinjlim_{\text{$L/K$ purely inseparable}} \quad \GL \right)^\wedge \\
&= \left( \varinjlim \, (\GK \stackrel{\sigma}{\to} \GK \stackrel{\sigma}{\to}
 \cdots ) \right)^\wedge \\
\Galg &= \left( \varinjlim_{\text{$L/K$ algebraic}} \quad \GL \right)^\wedge \\
&= (\Gsep \otimes_{\GK} \Gperf)^\wedge.
\end{align*}
From the canonicality of the constructions, we have the equalities of Galois groups
\begin{align*}
\Gal(L/K) &\cong
\Gal(\GL/\GK) \\
\Gal(K^{\sep}/K) &\cong
\Gal(\Gsep/\GK) \\
&\cong
\Gal(\Galg/\Gperf).
\end{align*}

Since $K^{\perf}$ and $K^{\alg}$ are perfect, their Witt rings
tensored with $\calO$ are isomorphic to $\Gperf$ and
$\Galg$, respectively, with $\sigma$ acting as the canonical
Frobenius on each of the Witt rings. 
To construct $\Gimm$, fix an
embedding of $K^{\alg}$ into $K^{\imm}$; such an embedding exists
because $K^{\imm}$ is algebraically closed. (See \cite{bib:me1} and
\cite{bib:me2} for a proof, and a description of the image of the
embedding.) Now define $\Gimm$ as the Witt ring of $K^{\imm}$
tensored with $\calO$ over $W$, and embed $\Galg$ into it
using the functoriality of the Witt vectors. Concretely, if one
chooses $u\in \Gimm$ in the image of the Teichm\"uller map
(that is, having $p^n$-th roots for all $n \in \NN$) and having residual
valuation 1,
one can
describe $\Gimm$ as the ring of series $\sum_{i \in \QQ} c_i
u^i$, with $c_i \in \calO$,
such that for each $\epsilon>0$, the set of $i$ with $|c_i| \geq
p^{-\epsilon}$ is well-ordered.

To extend Frobenius to $\Gimm$, first note that there exists
$u \in \Gperf$ such that $u^\sigma = u^p$; since
$\sigma$ is a bijection on $\Gperf$, $u$ is in the
image of the Teichm\"uller map. Now define
\[
\left(\sum_{i \in \QQ} c_i u^i\right)^\sigma =
\sum_{i \in \QQ} c_i^\sigma u^{pi}.
\]

In addition to Frobenius, we will also need to propagate derivations
on $\GK$ up to extensions, as far as is possible. Let $\theta$
be a derivation on $\GK$ over $\calO$.
For $L/K$
separable, $\theta$ extends uniquely to $\GL$, since the latter
can be written as $\GK[x]/(P(x))$ with $P(x)$ residually separable,
and we have $\theta(P(x)) = (\theta P)(x) + P'(x) \theta(x)$,
compelling us to set $\theta(x) = -(\theta P)(x) / P'(x)$. (Here $\theta P$
is the result of applying $\theta$ to the coefficients of $P$, while
$P'$ is the formal derivative.)

For $L/K$ inseparable, the situation is more complicated. For $c>0$ real,
we define $\Gamma^{\perf(c)}$ (resp.\ $\Gamma^{\alg(c)}$)
as the subring of $\Gperf$ (resp.\ $\Galg$)
consisting of those $x$ such that for each
$n \geq 0$, there exists $r_n \in \GK$ (resp.\ $\Gsep$)
such that $|x^{\sigma^{n}}-r_n|
< p^{-cn}$.
\begin{prop}
Any derivation $\theta$ on $\GK$ extends uniquely to a derivation
mapping $\Gamma^{\perf(c)}$ to $\Gamma^{\perf(c-1)}$
and $\Gamma^{\alg(c)}$ to $\Gamma^{\alg(c-1)}$, for each $c>1$.
\end{prop}
\begin{proof}
We will only give the arguments for $\Gamma^{\perf(c)}$, as the arguments
for $\Gamma^{\alg(c)}$ are the same.
Also, we
may assume $\theta$ does not map $\GK$ into $\pi\GK$.

Pick any residual uniformizer $t$ in $\GK$, and let $u =
\theta(t^\sigma) / \theta(t)^\sigma$. Then for all $x \in \GK$,
$\theta(x^\sigma) = u (\theta x)^\sigma$. (First check this for powers
of $t$, then extend by linearity.)

We claim that $|u| \geq p^{-1}$, which is to say $u \not\equiv 0 \pmod{\pi p}$.
Suppose that the contrary holds.
Recalling that $t^\sigma = t^p + \pi v$ for some $v \in \GK$,
we have $\theta(t^\sigma) = pt^{p-1} \theta(t) + \pi \theta(v)$. By hypothesis,
$\theta(x^\sigma) \equiv 0 \pmod{\pi p}$, which is to say
$\theta(v) \equiv - (p/\pi) t^{p-1}\theta(t) \pmod{p}$.

If we write
$v = \sum_i c_i t^i$, then $\theta(v) = \sum_i ic_i t^{i-1}\theta(t)$.
Since $\theta t$ is not divisible by $\pi$ by hypothesis, we must have
$ic_i \equiv 0 \pmod{p}$ for $i \neq p$ and $pc_p \equiv -p/\pi \pmod{p}$,
which is an absurdity. Thus $|u| \geq p^{-1}$ as claimed.

Now for $x \in \Gamma^{\perf(c)}$ and $n \geq 0$, define
$y_n \in \Gperf$ as follows. Choose $r_n \in \GK$
such that $|x^{\sigma^n} - r_n| < p^{-cn}$, and set
\[
y_n = \frac{\theta(r_n)^{\sigma^{-n}}}{u^{\sigma^{-1}} \cdots u^{\sigma^{-n}}}.
\]
For $n_1, n_2 \geq N$, we have $|y_{n_1} - y_{n_2}| < p^{-N(c-1)}$,
so the sequence $\{y_n\}$ converges in $\Gperf$ to a limit $y$
such that $|y - y_n| < p^{-n(c-1)}$.
Moreover,
\[
|y^{\sigma^n} - r_n u^{-n}| = |(y - y_n)^{\sigma^{n}}| < p^{-n(c-1)},
\]
so that $y \in \Gamma^{\perf(c-1)}$.
It is easily verified that the assignment $x \mapsto y$ yields a derivation,
that this derivation extends $\theta$, and that any derivation extending
$\theta$ must agree with this one on $\Gamma^{\perf(c)} \cap
(\GK)^{\sigma^{-n}}$ for each $n \in \NN$; from this final assertion
follows the uniqueness.
\end{proof}

\subsection{The ring $\Omega$ and the standard operators}
A number of our arguments will be simplified by using Frobenius and derivation 
operators of a particular simple form, called the \emph{standard operators};
we introduce these operators in this section. Some basic definition will have to
be made without this simplifying assumption, though, and the reasons for this will
become evident during the following constructions.
        
Let $t$ be a residual uniformizer in $\Gamma$. Then the subring of $\Gamma$
consisting of series $\sum_{n=0}^\infty c_n t^n$, with $c_n \in \calO$, will be
denoted $\Omega_t$; the decoration $t$ will be omitted if the choice of $t$ is to
be understood.

The Frobenius $\sigma_t$ mapping $t$ to $t^p$ will play a special role
in our work; we call it the \emph{standard Frobenius with respect to
$t$} (or simply the \emph{standard Frobenius} if $t$ is understood). We say
a Frobenius on $\Gamma$ is \emph{standard} if it equals $\sigma_t$ for some
residual uniformizer
$t \in \Gamma$. Similarly, we define the derivation $\theta_t = t \frac{d}{dt}$
mapping $t$ to itself and call it the \emph{standard derivation with respect to
$t$}. 
(Beware that in \cite{bib:dej1}, the standard derivation is $\frac{d}{dt}$; we
will comment further on the ramifications of this discrepancy.)

If the Frobenius $\sigma$ on $\Gamma$ is equal to $\sigma_t$,
we can regard $\Omega_t$ as a subring of $\Gamma$ stable under Frobenius.
This construction is not well-behaved under base change, however, which is to say that
given $L/K$, one cannot hope to construct $\Omega_u$ in $\GL$ and
$\Omega_t$ in $\GK$ such that $\Omega_u$ contains $\Omega_t$. The problem
is that the Frobenius given by $u^\sigma = u^p$ and the Frobenius given by
$t^\sigma = t^p$ never coincide, unless $L/K$ is a tamely ramified extension
followed by a purely inseparable extension, in which case one can take
$u = t^{1/n}$ for $n = [L:K]$. (In other words, a standard Frobenius
does not remain standard under a wildly ramified base extension.)

Note that over $\Gperf$, $\Galg$, or $\Gimm$, there must exist $u$
such that $u^\sigma = u^p$; in fact, this holds if and only if $u$ is
a Teichm\"uller element. Indeed, if $u^{\sigma} = u^p$, then
$u^{\sigma^{-n}}$ is a $p^n$-th root of $u$ for each $n \in \NN$, so
$u$ must be a Teichm\"uller element. Conversely, if $u$ is a
Teichm\"uller element, it maps to a Teichm\"uller element under every ring
endomorphism, in particular under $\sigma$; but $u^\sigma \equiv u^p
\pmod{\pi}$ and both sides are Teichm\"uller elements, so we must have equality. In short, in a situation where we work only over a ring containing
$\Gperf$, we may assume every Frobenius is standard. (In particular, one
can canonically define $\Omega^{\perf}, \Omega^{\alg}, \Omega^{\imm}$,
as the Witt rings of the valuation rings within their respective residue
fields, but we will not refer to these rings elsewhere.)

\subsection{Overconvergent rings}
Our next task is to construct ``overconvergent'' subrings of the rings
$\Gamma^{K,*}$ defined in the previous section, and to determine to what
extent the Frobenius and derivation operators extend to these subrings. In particular,
given a Frobenius on $\GK$ which is overconvergent with respect to some residual
uniformizer,
we will define a subring $\Gamma^*_{\con}$ of $\Gamma^{*}$ for $* \in \{\emptyset, L,
\sep, \perf, \alg, \imm\}$ mapped into itself by the extension of Frobenius to
$\Gamma^*$. Moreover, these subrings will be compatible in that if $*_1$ and $*_2$
are decorations such that $\Gamma^{*_1} \subseteq \Gamma^{*_2}$,
then $\Gamma^{*_1}_{\con} \subseteq \Gamma^{*_2}_{\con}$.

We first construct overconvergent subrings of $\Gamma$.
If $t$ is a residual uniformizer in $\Gamma$,
we can construct the functions
$v_{t,n}$ for $n \geq 0$
by expressing any $u \in \Gamma$ as $u=\sum_i u_i t^i$ and
letting $v_{t,n}(u)$ be the smallest integer $j$ such that $|u_j| \geq p^{-n}$,
or $\infty$ if $j$ does not exist.
(We do not require $n$ to be integral, but we may as well force it to be an
integral multiple of the integer $m$ with the property that $|\pi| = p^{-1/m}$.)

\begin{prop} \label{prop:vtn}
The following relations are satisfied by the functions $v_{t,n}$ on $\Gamma$,
and determine them uniquely.
\begin{enumerate}
\item For $n \geq 0$ and $x,y \in \Gamma$, $v_{t,n}(x+y) \geq \min\{v_{t,n}(x), v_{t,n}(y)\}$,
with equality if $v_{t,n}(x) \neq v_{t,n}(y)$.
\item For $n \geq 0$ and $x,y \in \Gamma$, $v_{t,n}(xy) \geq \min_{i\geq 0} \{v_i(x) + v_{n-i}(y)\}$.
\item
For any residual uniformizer $t$,
$v_{t,0}(x)$ equals the valuation of the residue of $x$
in $K$.
\item
For $n \geq 0$, $v_{t,n}(t) = 1$.
\item
For $n \geq 0$ and $x \in \Gamma$, $v_{t,n}(x^{\sigma_t}) = pv_{t,n}(x)$.
\end{enumerate}
\end{prop}
We say that $u$ is \emph{overconvergent} with respect to a residual
uniformizer $t$
if $v_{t,n}(u) \geq -cn-d$ for some constants $c,d > 0$;
this relation is symmetric if $u$ is also a residual uniformizer.
We say a Frobenius $\sigma$ on $\Gamma$ is overconvergent with respect to
$t$ if $t^\sigma$ is overconvergent; likewise, we say a derivation $\theta$
on $\Gamma$ is overconvergent with respect to $t$ if $\theta(t)$ is overconvergent.
In particular, the standard operators $\sigma_t$ and $\theta_t$ are overconvergent
with respect to $t$.
It will follow from Proposition~\ref{prop:frobcon} (see below) that
if $t_1$ and $t_2$
are residual uniformizers with respect to which $\sigma$ is overconvergent,
then $t_1$ and $t_2$ are overconvergent with respect to each other.

Let $\GK_{\con(t)}$ be the subring of $\GK$ consisting of those elements
of $\GK$ which are overconvergent with respect to $t$; we will drop $t$
from the notation when it is understood. Elsewhere in the literature, this
subring is notated using a dagger; this is a throwback to Monsky-Washnitzer
dagger cohomology, whence the notion of overconvergence originated.

\begin{prop}
Suppose that in Proposition~\ref{prop:extgam}, $R = \Gamma$ and $\sigma$ is 
overconvergent with respect to some residual uniformizer $t$. Then
$\sigma$ on $S$ is overconvergent with respect to some residual uniformizer
$u$, and (viewed as an element of $S$) $t$ is overconvergent with respect to $u$.
\end{prop}
\begin{proof}
It suffices to consider the cases in which $L/K$ is separable but tamely ramified,
purely inseparable, or an Artin-Schreier extension, as any $L/K$ can be expressed
as a tower of these. In the first case, we can choose $u = t^{1/n}$, where $n = [L:K]$,
and $u^\sigma = (t^{\sigma})^{1/n}$ is clearly overconvergent with respect to $u$.
In the second case, we can choose $u$ such that $u^\sigma = t$.

In the third case, choose $u$ such that $u^{-p} - u^{-1} = t^{-k}$ for
some $k \in \NN$ not divisible by $p$; then $u^{1/k}$ exists and is a
residual uniformizer, with respect to which $t$ is overconvergent. Now
$u^\sigma$ is a root of the polynomial $P(x) = x^p + x^{p-1}
(t^{\sigma})^k + (t^{\sigma})^k$, which has coefficients which are
overconvergent with respect to $u$. However, Crew
\cite[Proposition~4.2]{bib:crew1} has shown that the ring of
overconvergent series with respect to a given residual uniformizer is
henselian. Thus $u^\sigma$ is also overconvergent with respect to $u$.
\end{proof}
In symbols, this says there exists $u$ such that $\GL_{\con(u)}$ contains
$\GK_{\con(t)}$; when $t$ is understood, we will suppress $u$ as well and simply
say that $\GL_{\con}$ contains $\GK_{\con}$.

We now wish to extend the notion of overconvergence to $\Gsep$,
$\Gperf$, $\Galg$ and $\Gimm$, assuming that overconvergence
is taken with respect
to some $t \in \Gamma$ with respect to which $\sigma$ is overconvergent.
Unfortunately, it is not obvious how to extend $v_{t,n}$ to these larger
rings, except when $\sigma = \sigma_t$,
so we must use a somewhat indirect approach. (Note that this discussion is elided
in \cite{bib:dej1}; specifically,
the existence of $v_{t,n}$ for $\sigma = \sigma_t$ is stated and used
but not justified.)

Recall from the previous section that any Teichm\"uller element $u \in \Gperf$
has the property that $u^\sigma = u^p$. Let $u$ be the Teichm\"uller lift
of an element $K^\perf$ of valuation 1, and define the function $v_{u,n}$
on $\Gimm$ mapping $x = \sum c_i u^i$
to the smallest rational number $j$ such that
$|c_j| \geq p^{-n}$; then the analogue of Proposition~\ref{prop:vtn}
holds with $v_{t,n}$ replaced by $v_{u,n}$.
Define $\Gimmcon$ as the set of series $x = \sum c_i u^i$
such that $v_{u,n}(x) \geq -cn-d$ for some constants $c,d>0$, and set
$\Gsepcon = \Gsep \cap \Gimmcon$ and so forth. By construction, $\sigma$
maps $\Gamma^{*}_\con$ into itself.

There are a number of compatibilities that must be verified for the
above definition. In particular, we must show that $v_{u,n}$, and hence
the overconvergent rings, depend neither on the choice of
the Teichm\"uller element $u \in \Gperf$, nor on the embedding of
$K^{\alg}$ into $K^{\imm}$.
Fortunately, we can give an alternate characterization of $v_{u,n}$
on $\Gimm$ from which both of these compatibilities are manifest.
\begin{prop}
Let $\tau$ denote the Teichm\"uller map from $K^{\imm}$ to $\Gimm$,
and let $\nu$ denote the valuation on $K^{\imm}$, normalized so that
a uniformizer of $K$ has valuation $1$.
Also let $|\pi| = p^{-1/m}$.
If $x = \sum_{n=0}^\infty \pi^n \tau(d_n)$ with $d_n \in K^{\imm}$,
then
\[
v_{u,-n/m}(x) = \min \{ \nu(d_0), \dots, \nu(d_{n-1}) \}.
\]
\end{prop}
\begin{proof}
We first establish that $v_{u,n/m}(\tau(x)) = \nu(x)$ for all $x \in K^{\imm}$.
From the definition of the Teichm\"uller map, $\tau(x)$ is congruent
modulo $\pi^n$ to $(y^{\sigma^{-n}})^{p^n}$ for any $y \in \Gimm$ lifting $x$.
In particular, we can pick $y$ to be a series $\sum c_i u^i$ such that
$c_i = 0$ for $i<\nu(x)$, and then $(y^{\sigma^{-n}})^{p^n}$ will have the same property,
so $v_{u,n/m}(\tau(x))$ will not be less than $\nu(x)$ (and not greater, since
$v_{u,0}(\tau(x)) = \nu(x)$).

To establish the desired formula for $x$ arbitrary, let $j = \min\{\nu(d_0), \dots,
\nu(d_{n-1})\}$, and let $k$ be the smallest nonnegative integer such that
$\nu(d_k) = j$. On one hand, we
have $v_{u,-n/m}(x) \geq j$ by the previous paragraph together with
Proposition~\ref{prop:vtn} (or rather,
by its analogue for $v_{u,n}$). On the other hand, modulo $\pi^n$,
the coefficient of $u^j$ receives zero contribution from $\tau(d_0), \dots,
\pi^{k-1} \tau(d_{k-1})$, a contribution divisible by $\pi^k$ but not by
$\pi^{k+1}$ from $\tau(d_k)$, and contributions divisible by $\pi^{k+1}$
from $\tau(d_{k+1}), \dots, \tau(d_{n-1})$. Therefore the coefficient of $u^j$
in $x$ is nonzero modulo $\pi^n$, and $v_{u,-n/m}(x) \leq j$. We conclude
$v_{u,-n/m}(x) = j$, proving the desired result.
\end{proof}
Since the formula on the right hand side does not involve $u$, the independence
of $v_{u,n}$ from $u$ is immediate. As for the independence from the embedding
of $\Galg$ into $\Gimm$, recall that
the valuation $\nu$
on $K^{\alg}$ is Galois-equivariant (a standard fact about local
fields; see \cite[Chapter~I]{bib:serre} for proof), as is the decomposition of
an arbitrary element into Teichm\"uller elements.
Since we now have that $v_{u,n}$ is Galois-invariant
on $\Galg$, we deduce the Galois-invariance of
$\Galgcon$ and $\Gsepcon = \Gsep \cap \Galgcon$ as well.

Finally, we must show that $\GL \cap \Gimmcon = \GL_{\con}$;
it suffices to establish this for $L=K$, which follows from the following proposition.
(More precisely, this proposition establishes that $\Gcon \subseteq \Gimmcon$;
the argument that $x \in \Gamma \cap \Gimmcon$ implies $x \in \Gcon$ is similar.)
\begin{prop} \label{prop:frobcon}
Suppose the Frobenius $\sigma$ on $\Gamma$
is overconvergent with respect to some residual uniformizer
$t \in \Gamma$. Then for any $u \in \Gperf$ with residual valuation $1$ such that $u^\sigma
= u^p$, we have $v_{u,n}(t) \geq -cn$ for suitable $c>0$.
\end{prop}
\begin{proof}
Let $e>0$ be such that $v_{t,n}(t^\sigma) \geq -en$ for $n \geq 0$. (This is possible
because $\sigma$ is overconvergent with respect to $t$ and because
$v_{t,n}(t^\sigma) = p > 0$.)
We will construct $t_m \in \Gperf$ for $m \geq 0$ such that
$v_{u,n}(t_m) \geq -en/p$ for $m,n \geq 0$, and $t_m \equiv t \pmod{\pi^m}$; the existence
of such $t_m$ suffices to prove the desired assertion with $c=e/p$.

We may start with $t_0 = u$. Now suppose $t_0, \dots, t_m$ have been constructed.
Write $t^{\sigma} = t^p + \pi \sum_{i \in \ZZ} c_i t^i$, and set
\[
t_{m+1} = \left( t_m^p + \pi \sum_{i \in \ZZ} c_i t_m^i \right)^{\sigma^{-1}}.
\]
Then
\begin{align*}
t_{m+1}^{\sigma} &= t_m^p + \pi \sum_{i \in \ZZ} c_i t_m^i \\
&\equiv t^p + \pi \sum_{i \in \ZZ} c_i t^i = t^\sigma \pmod{\pi^{m+1}},
\end{align*}
so that $t_{m+1} \equiv t \pmod{\pi^{m+1}}$. On the other hand,
$v_{u,n}(t_m^p) \geq -en/p$, and
\begin{align*}
v_{u,n}(c_i t_m^i) &= v_{u,n}(c_i u^i (t_m/u)^i) \\
&\geq \min_{j \leq n} \{v_{u,j}(c_i u^i) + v_{u,n-j} (t_m/u)^i \} \\
&\geq \min_{j \leq n} \{-ej - e(n-j)/p\} = -en
\end{align*}
(the estimate $v_{u,j}(c_i u^i) \geq -ej$ following from the bound $v_{t,j}(t^\sigma) \geq -ej$).
Thus $v_{u,n}(t_{m+1}) \geq en/p$ for all $n \geq 0$, and the construction of $t_{m+1}$
is complete.
\end{proof}

The following basic lemma is an extension of \cite[Proposition
8.1]{bib:dej1}.
\begin{prop} \label{inj1}
For $\sigma = \sigma_t$ standard, the following multiplication maps are injective:
\begin{eqnarray*}
\Gperfcon \otimes_{\Gcon} \Gamma &\to& \Gperf \\
\Gsepcon \otimes_{\Gcon} \Gamma &\to& \Gsep \\
\Galgcon \otimes_{\Gcon} \Gamma &\to& \Galg \\
\Gimmcon \otimes_{\Gcon} \Gamma &\to& \Gimm \\
\Gperfcon \otimes_{\Gcon} \Gsep &\to& \Galg.
\end{eqnarray*}
\end{prop}
\begin{proof}
The first, second, and third assertions
follow from the fourth one, since in the diagram
\[
\xymatrix{
\Gamma^*_{\con} \otimes_{\Gcon} \Gamma \ar[r] \ar[d] & \Gamma^{*} \ar[d] \\
\Gimmcon \otimes_{\Gcon} \Gamma \ar[r] & \Gimm
}
\]
(with $* \in \{\perf, \sep, \alg\}$)
the left vertical arrow is injective by flatness ($\Gcon \to
\Gamma$ is an unramified extension of discrete valuation rings,
hence flat). Thus we concentrate our attention on proving that
$\Gimmcon \otimes_{\Gcon} \Gamma \to \Gimm$ is injective.

Suppose $\sum_{i=1}^n f_i \otimes g_i$ is a nonzero element of
$\Gimmcon \otimes_{\Gcon} \Gamma$ such that $\sum f_i g_i = 0$ in
$\Gimm$, and such that $n$ is minimal for the existence of such an element.
Then the $g_i$ are linearly independent over $\Gcon$, otherwise
we could replace one of them by a combination of the others and
decrease $n$.

Now each $x \in \Gimm$ can be uniquely written in the form
$\sum_{0 \leq \alpha < 1} x_{\alpha} t^{\alpha}$,
and $x \in \Gimmcon$ implies $x_{\alpha} \in \Gcon$
for all $\alpha$. Writing
$f_i = \sum f_{i,\alpha} t^{\alpha}$, we find (by the uniqueness
of the decomposition of $\sum f_i g_i = 0$) that
$\sum_{i=1}^n f_{i,\alpha} g_i = 0$ for each $\alpha$. However,
since the $g_i$ are linearly independent over $\Gcon$, we have
$f_{i,\alpha} = 0$ for all $i$ and $\alpha$, contradicting
the fact that $\sum f_i \otimes g_i$ is nonzero.

To prove the final assertion, suppose $\sum_{i=1}^n f_i \otimes
g_i$ is a nonzero element of $\Gperfcon \otimes_{\Gcon} \Gsep$
such that $\sum f_ig_i = 0$ in $\Galg$, and such that $n$ is
minimal for the existence of such an element.  By the first
assertion, we have that $\sum_{i=1}^n f_i \otimes g_i$ maps to
zero in $\Gperf \otimes_{\Gamma} \Gsep$, which means that the
$f_i$ are linearly dependent over $\Gamma$. Choose $r_1, \dots,
r_n$ in $\Gamma$, not all zero, such that $\sum r_i f_i = 0$.
Without loss of generality, suppose that $r_1 \neq 0$
and $|r_1| \geq |r_i|$ for all $i$.  Then
$\sum_{i=2}^n f_i \otimes (g_i - g_1 r_i/r_1)$ maps to zero in
$\Galg$, and by the minimality of $n$, we must have $g_i = g_1
r_i/r_1$ for all $i$, which is to say $\sum f_i \otimes g_i = (1
\otimes g_1/r_1) \sum f_i \otimes r_i$. Since the $r_i$ lie in
$\Gamma$, we can apply the second assertion to deduce that $\sum
f_i \otimes r_i = 0$ in $\Gperfcon \otimes_{\Gcon} \Gamma$ and
hence also in $\Gperfcon \otimes_{\Gcon} \Gsep$.  Thus $\sum f_i
\otimes g_i=0$ as well.
\end{proof}

We will often use the fact that Galois descent works for overconvergent rings,
so let us state this explicitly.
\begin{prop} \label{prop:galdesc}
Let $M$ be a finitely generated free module over $\Gcon$ (resp.\ $\Gperfcon$)
and $X$ a Galois-stable submodule of $M \otimes_{\Gcon} \Gsepcon$
(resp.\ $M \otimes_{\Gperfcon} \Galgcon$). Then $X$ can be expressed as
$Y \otimes_{\Gcon} \Gsepcon$ (resp.\ $Y \otimes_{\Gperfcon} \Galgcon$) for
some submodule $Y$ of $M$.
\end{prop}
\begin{proof}
Suppose $\be_1, \dots, \be_n$ is a basis of $M$.
After reordering the $\be_i$ suitably, one can find a basis
$\bv_1, \dots, \bv_m$ of $X$ such that if one writes
$\bv_i = \sum_j c_{ij} \be_j$ with $c_{ij} \in \Gsepcon$,
then $c_{ij}=0$ if $j \leq m$ and $i \neq j$ (by Gaussian elimination).
Now simply replace $\bv_i$ by $\bv_i/c_{ii}$ and one gets a set of 
vectors in $M$, and we can take $Y$ to be their span.
\end{proof}
Another formulation, which we will often invoke, is that $M$ is a finitely
generated free module over $\Gcon$ and an element of
$\wedge^k M$ factors completely over $\Gsepcon$, then said element factors
completely over $\Gcon$.

Before concluding the discussion of overconvergent rings, one caveat
must be made about derivations on overconvergent rings. While a
derivation on $\Gamma$ extends to $\Gsep$, a derivation on $\Gamma$
that carries $\Gcon$ into itself need not do likewise on $\Gsepcon$.
On the other hand, for $c>1$, the extension of a derivation on
$\Gamma$ to a derivation from $\Gamma^{\perf(c)}$ to
$\Gamma^{\perf(c-1)}$ does map $\Gamma^{\perf(c)}_{\con}$ into
$\Gamma^{\perf(c-1)}_{\con}$.

\subsection{Analytic rings}

In this section, we introduce some rings with the decoration ``$\an$'', which
mostly correspond to rings of rigid analytic functions on certain regions
of a $p$-adic analytic space. Unlike the other rings introduced so far,
these rings are not discrete valuations rings; in fact, they are not
even local.

The ring $\Oan$ consists of those series $\sum_{n=0}^{\infty} c_{n}
t^{n}$ in $\calO[\fp][[t]]$ such that $\limsup_{n \to \infty}
|c_n|^{1/n} \leq 1$.
These series can be identified with the rigid analytic
functions on the formal unit disc. Given a residual uniformizer $t$ in
$\GK$,
the set $\GK_{\an}$ is analogously defined
as the set of series $\sum_{n \in \ZZ}c_n t^n$ such that
$\limsup_{n \to \infty} |c_n|^{1/n} \leq 1$ and
$|c_n| \to 0$ as $n \to -\infty$. However, $\Gan$ cannot be made into a ring
using series addition and multiplication: attempting to multiply two series
in $\Gan$ can lead to expressions for the coefficients of the result involving
infinite sums. Fortunately, one can define the subset $\Gancon$ of $\Gan$
consisting of series $\sum c_n t^n$ such that $|c_n| \leq p^{-cn-d}$ for 
some constants $c,d>0$ (depending on the series), and $\Gancon$ is a ring.
This definition depends only on the choice of $\Gcon$ with $\Gamma$ and
not on $t$, in the sense that
$\Gancon$ remains unchanged if it is defined in terms of another
uniformizer $t_1$ which is overconvergent with respect to $t$.
Also, $\Gan$ is naturally a module over $\Gcon$ (though not over $\Gancon$ or
$\Gamma$).

We will have occasion to extend the functions $v_{t,n}$, originally
defined on $\Gamma$, to $\Gan$ and $\Gancon$. In fact, we will use
$v_{t,n}$ for $n$ arbitrary (not necessarily nonnegative),
defined as mapping $x = \sum c_n t^n \in \Gan$ to the smallest integer
$j$ such that $|c_j| \geq p^{-n}$, or $\infty$ if no such $j$ exists.
These maps satisfy a slightly modified version of Proposition~\ref{prop:vtn},
given below.
\begin{prop} \label{prop:vtn2}
The following relations are satisfied by the functions $v_{t,n}$
on $\Gan$.
\begin{enumerate}
\item
For $n \geq 0$ and $x \in \Gamma$, $v_{t,n}(x)$ takes the same value
whether evaluated in $\Gamma$ or in $\Gan$.
\item For $n \in \QQ$ and $x,y \in \Gan$, $v_{t,n}(x+y) \geq \min\{v_{t,n}(x), v_{t,n}(y)\}$,
with equality if $v_{t,n}(x) \neq v_{t,n}(y)$.
\item For $n \in \QQ$, $x \in \Gan$ and $y \in \Gancon$, $v_{t,n}(xy) \geq \min_{i \in \QQ} \{v_i(x) + v_{n-i}(y)\}$\
.
\item
For $n \in \QQ$, $v_{t,n}(t) = 1$.
\item
For $n \in \QQ$ and $x \in \Gan$, $v_{t,n}(x^{\sigma_t}) = pv_{t,n}(x)$.
\end{enumerate}
\end{prop}

For a finite extension $L$ of $K$, it is easy to see that the
embedding of $\GK_{\con}$ into $\GL_{\con}$ extends
canonically to an embedding of $\GK_{\an,\con}$ into
$\GL_{\an,\con}$, and that automorphisms of $L$ over $K$ give rise
to automorphisms of $\GL_{\an,\con}$ over $\GK_{\an,\con}$.
Dealing with larger extensions is more
complicated; as we did for convergent rings, we begin at the top
with $\Gimm$. For $u \in \Gperf$ a Teichm\"uller element with residual
valuation 1,
the ring $\Gimm_{\an,\con}$ is defined as the set of generalized power
series $\sum_{i \in \QQ} x_i u^i$, with $x_i \in \calO[\fp]$, such
that for each $\epsilon \in \RR$, the set of $i$ such that $|x_i| \geq
p^{\epsilon}$ is well-ordered.
As in the definition of $\Gimmcon$, the choice of $u$ does not affect
the definition.

We define $\Gamma^{*}_{\an,\con}$ for $* \in \{\sep, \perf, \alg\}$
as the set of $x \in \Gimmcon$ such that for every $m,n \in \QQ$,
$x$ can be written as an element of $\Gamma^{*}[\fp]$, plus
$\pi^m$ times an element of $\Gimm$, plus a series of the form
$\sum_{i \geq n} x_i u^i$. One can verify that this indeed gives
a subring, and that $\Gal(K^{\alg}/K^{\perf})$ acts on $\Galg_{\an,\con}$
and $\Gsep_{\an,\con}$. Beware, though,
that $\Gamma^{*}_{\an,\con} \cap \Gimm$ is not equal to $\Gamma^{*}_{\con}$;
it is actually a larger ring with residue field equal to the completion of
$K^{*}$ (in the valuation induced from $K$).
In particular, the completions of $K^{\sep}$ and $K^{\alg}$ are equal, so
$\Gsep_{\an,\con} = \Galg_{\an,\con}$.


One convenient feature of $\Gancon$, not shared by any of the other
rings introduced so far, is that not only does a derivation on $\Gcon$
extend to $\Gancon$, but said derivation admits an antiderivative as
well. More precisely, for $x \in \Gancon$, there exists $y \in
\Gancon$ such that $\theta_t y = x$ if and only if the constant
coefficient of $x$ is zero. (An analogous statement for a non-standard
derivation can be made using the fact that any such derivation is a
scalar multiple of $\theta_t$.)

\section{Crystals and their properties}

In this section we define crystals, as needed for our purposes, and describe
their basic structural properties. Our format follows Katz \cite{bib:katz}
and de~Jong \cite{bib:dej1} with some minor modifications. Note that we use
the term ``crystal'' where other sources use ``isocrystal''; as we work
entirely in the local setting, this should not cause any ambiguity. 

\subsection{Crystals}

Let $R$ be any characteristic 0
ring from the previous section, and $\sigma$ a Frobenius on $R$.
An \emph{$F$-crystal} over $R$ is a 
finite, locally free $R$-module $M$ equipped with an additive, $\sigma$-linear
endomorphism $\map FMM$ which becomes an isomorphism over $R[\fp]$.
More precisely, if we put $M^\sigma = M \otimes_{R,\sigma} R$, $F$
should be an $R$-linear map $\map{F}{M^\sigma}{M}$ which becomes an isomorphism
after tensoring (over $\calO$) with $\calO[\fp]$.
Even more precisely, there should exist $\ell \in \NN$ such that
$p^\ell$ annihilates the kernel and cokernel of $F$.
The fundamental examples are the
\emph{trivial $F$-crystals}, which are rank one modules of the form $M=R$
with $F(x) = cx^\sigma$ for some $c \in \calO$.

Note that some sources (like~\cite{bib:tsu3}) allow crystals in
which $F$ is a $\sigma^k$-linear endomorphism for $k \geq 1$. For our
purposes, it suffices to note that one can make such an object into a
crystal in our sense at the expense of multiplying its dimension by $k$.
Namely, given the $R$-linear map $\map {F}{M^{\sigma^k}}{M}$, then
\begin{align*}
F_1: M^{\sigma^{k}} \oplus M^{\sigma^{k-1}} \oplus \cdots \oplus M^{\sigma}
&\to
M^{\sigma^{k-1}} \oplus \cdots \oplus M^\sigma \oplus M \\
(m_k, \dots, m_1) &\mapsto (m_{k-1}, \dots, m_1, F(m_k))
\end{align*}
is $R$-linear as well.

Suppose $R$ is a ring admitting a derivation $\theta$ over $\calO$.
Then we say $M$ is an
\emph{$(F, \nabla)$-crystal} over $R$
if $M$ is an $F$-crystal over $R$
equipped with an $\calO$-linear connection $\nabla$, that is,
an additive map $\map{\nabla}{M}{M}$ with the following properties:
\begin{enumerate}
\item (Leibniz rule) For all $a \in R, m \in M$,
$\nabla(am) = \theta(a)m + a \nabla(m)$.
\item (Frobenius compatibility)
For all $m \in M$, $\theta(t^\sigma)
F\nabla (m) = \theta(t)^\sigma \nabla F(m)$.
\end{enumerate}
In particular, the trivial $F$-crystals are also $(F,\nabla)$-crystals
using the connection $\nabla(x) = \theta(x)$.

While the notion of an $F$-crystal depends strongly on the choice of
$\sigma$ (though we will see how to get around this choice below),
over $\Gamma$ the choice of $\theta$ is comparatively immaterial. Specifically,
if $\theta_1$ is another derivation, and assuming neither $\theta$ nor
$\theta_1$ maps $\Gamma$ into $\pi \Gamma$,
then there exists $c$ such that
$\theta_1(x) = c\theta(x)$ for all $x \in R$, namely
$c = \theta_1(t)/\theta(t)$. Now setting $\nabla_1(x) = c \nabla(x)$
gives a new connection satisfying the revised Leibniz rule
$\nabla_1(am) = \theta_1(a)m + a \nabla_1(m)$.

On the other hand, over $\Omega$ the choice of $\theta$ in the definition
of $(F, \nabla)$-crystal is quite significant. We will ordinarily use the
standard derivation $\theta_t$ mapping $t$ to $t$; this is a departure from
\cite{bib:dej1}, in which the derivation $\frac{d}{dt}$ mapping $t$ to 1 is
used. The result is that we allow as crystals certain objects that originate
in geometry as ``crystals with logarithmic poles''; this permissiveness will
be crucial for the correct statement of the semistable reduction conjecture.

The study of crystals in many ways resembles a ``$\sigma$-twisted'' analogue of
ordinary linear algebra, and some of our terminology will reflect this
resemblance.
For example, a nonzero element $\bv$ of an $F$-crystal $M$ over $R$ is said to be an
\emph{eigenvector} if there exists $\lambda \in \calO$ such that $F\bv =
\lambda \bv$.  The $p$-adic valuation of $\lambda$ is called the
\emph{slope} of the eigenvector. (Note that a scalar multiple of an
eigenvector is not ordinarily an eigenvector, unless the scalar lies in $\calO$.)

A \emph{morphism} between $F$-crystals (resp.\ $(F, \nabla)$-crystals) $M_1$ and $M_2$ over $R$ is
an $R$-linear map from $M_1$ to $M_2$ which makes the obvious diagrams commute.
A morphism $f: M_1 \to M_2$ is said to be an \emph{isomorphism} (resp. \emph{isogeny})
if there exists a morphism $g: M_2 \to M_1$ such that $f \circ g$ and $g \circ f$
are the identity maps (resp.\ are the same scalar
multiple of the identity maps) on their respective domains. We
will always work in the category of crystals up to isogeny, which
is to say what we call a ``crystal'' is in reality an isocrystal.
Since we will never consider the category of crystals up to
isomorphism, we have lightened the notational load by dropping the
prefix ``iso'' throughout.

We will frequently encounter sets of elements of $M$ which form a basis
for $M \otimes_{\calO} \calO[\fp]$; such sets will be called
\emph{isobases} of $M$. Given an isobasis of $M$, we will refer frequently to
the matrices through which $F$ and (if applicable) $\nabla$ act on the isobasis.
Beware that \emph{a priori} these matrices only have entries in $R[\fp]$.
On the other hand, we will call two isobases \emph{commensurate} if they
generate the same submodule of $M$ over $R$; then if $F$ or $\nabla$ acts
on an isobasis through an integral matrix, obviously it acts on any
commensurate isobasis through another integral matrix.

Several other standard constructions of linear algebra carry over to crystals
without difficulty, such as tensor products and exterior powers.
Subcrystals are defined in the obvious manner; quotients objects are
defined by modding by $\pi$-power torsion in the naive quotient.
Duals require a bit more care:
to give $M^* = \Hom(M, R)$ the structure of an $F$-crystal, recall
that there exists $\ell$ such that $p^\ell$ annihilates the kernel and
cokernel of $F$. That means there exists $\map{\phi}{M}{M^\sigma}$ such
that $F \circ \phi$ and $\phi \circ F$ act by multiplication by $p^{2\ell}$.
Now the transpose of $\phi$ maps $(M^*)^\sigma$ to $M^*$, giving $M^*$ the
structure of an $F$-crystal. Beware that this structure is only well-defined
up to tensoring with a trivial crystal, because of the freedom in choosing
$\ell$; in case it is necessary to recall $\ell$, one may notate the dual
as $M^*(\ell)$.

The study of crystals is simplest when $R$ is $p$-adically complete
and $R/pR$ is an algebraically closed
field, thanks to the classification of Dieudonn\'e-Manin. (See \cite{bib:katz}
for details.)
Over such a ring $R$, after making a suitable totally ramified extension,
every $F$-crystal admits an isobasis of of eigenvectors; the slopes of these
eigenvectors are
called the slopes of the crystal.  The \emph{Newton polygon} of a
crystal of rank $n$ is the graph of the piecewise linear function from
$[0,n]$ to $\RR$ sending 0 to 0, whose slope between $k-1$ and $k$
equals the $k$-th smallest slope of the crystal (counting
multiplicities). This polygon turns out to be an isogeny invariant of the
crystal. (Katz also associates a second set of slopes to a
crystal, its \emph{Hodge slopes}. As these are not isogeny-invariant, we will not
discuss them here.)
We may extend the definitions of slopes and Newton polygons to a crystal over
any ring whose residue ring is a field, by extending scalars to obtain an
algebraically closed residue field. A crystal over $\Omega$ has two sets of
slopes: its \emph{generic slopes}, obtained by changing base to $\Gamma$,
and its \emph{special slopes}, obtained by changing base to $\calO$ by
reduction modulo $t$. (The special Newton polygon never goes below the
generic Newton polygon, by Grothendieck's specialization theorem
\cite[Theorem~2.3.1]{bib:katz}.)
On the other hand, a crystal over $\Gancon$ cannot
be given a meaningful set of slopes, because the base ring has no $p$-adic
valuation.

An $F$-crystal (resp.\ $(F, \nabla)$-crystal)
$M$ is \emph{unipotent} if after a suitable extension of
$\calO$, it becomes isogenous to a crystal admitting a filtration
$0 = M_0 \subset M_1 \subset \cdots \subset M_n$ of sub-$F$-crystals
(resp.\ sub-$(F,\nabla)$-crystals),
such that each successive quotient is trivial.
A crystal over $R$ is \emph{quasi-unipotent} if it becomes unipotent
after a finite separable extension of $R$.
An $F$-crystal $M$ is \emph{constant} if after a suitable extension of $\calO$,
it becomes isogenous to a direct sum of trivial crystals.
A constant crystal is unipotent, but not vice versa.

We say a crystal is \emph{isoclinic of slope $i$} if all of its slopes
are equal to $i$. (For $R = \Omega$, this includes both the special and
generic slopes, but the specialization theorem ensures that if one set
of slopes are all equal, so are the other set.)
For the rings considered in this paper, a result of Katz
\cite[Theorem~2.6.1]{bib:katz} (see also \cite[Lemma~6.1]{bib:dej1})
implies that for $\lambda \in \calO$,
a crystal is isoclinic of slope $-(\log_p |\lambda|)$ if and
only if it is isogenous to a crystal on which the action of $F$
factors through multiplication by $\lambda$.
Over $\Omega$, every isoclinic crystal 
is constant: assuming the slopes are all 0, choose
linearly independent elements on which $F$ acts by an integral matrix
$A$ congruent to 1 modulo $t$. (The congruence modulo $t$ uses the
Dieudonn\'e-Manin classification, and in particular relies on the
residue field being algebraically closed.) Then the infinite product $U =
 A A^\sigma A^{\sigma^2}\cdots$ converges $t$-adically and
$AU^\sigma = U$, so changing basis by $U$ gives a basis on which $F$ acts
through the identity matrix.

\subsection{Change of Frobenius}

The category of crystals over $R$ ostensibly depends on the choice of
a Frobenius. It is a fundamental property of $(F, \nabla)$-crystals 
(and one which is natural from the geometric perspective) that
this dependence is actually illusory; the proof below is due to
Tsuzuki \cite[Theorem~3.4.10]{bib:tsu3},
but we have shored up the justification
of a key point in the original argument.
(It is not known whether the same result holds for $F$-crystals; this
deficiency forces us to assume certain objects are $(F, \nabla)$-crystals
when other considerations only demand that they be $F$-crystals.)
\begin{prop} \label{prop:changefrob}
For $R$ equal to one of $\Gamma$ or $\Gcon$,
let $\sigma_1$ and $\sigma_2$ be two choices of Frobenius on $R$. Then
the categories of $(F, \nabla)$-crystals over $R$ equipped with
$\sigma_1$ and over $R$ equipped with $\sigma_2$ are equivalent.
\end{prop}

The proof of this assertion requires a technical lemma, regarding the
convergence of a certain sequence formed from $\nabla$. The reader is
advised to skip ahead to the proof of the proposition on first reading.
\begin{lemma} \label{lem:changefrob2}
Let $M$ be an $(F, \nabla)$-crystal over $R$. Then for every $\epsilon>0$,
there exists $c$ such that for all $\bv \in M$,
\[
\left| \frac{1}{n!} \nabla(\nabla - 1)\cdots (\nabla-n+1)\bv\right| < p^{n\epsilon+c}.
\]
\end{lemma}
\begin{proof}
For short, we write $\nabla^{(n)}$ for $\nabla(\nabla - 1)\cdots (\nabla-n+1)$
(with $\nabla^{(0)} = 1$).
Choose $\ell$ such that $p^{\ell}(p-1)\epsilon > 1$,
and
choose linearly independent elements of
$M$ on which $\nabla$ acts via a matrix $N$ with
$|N-1| < |p^{\ell}!|$.
This step is accomplished using $F$: given
$\be_1, \dots, \be_n$ on which $\nabla$ acts via a matrix $N$, $\nabla$
acts on $F\be_1, \dots, F\be_n$ via the matrix $N^\sigma
\theta(t^\sigma)/\theta(t)^\sigma$. Thus repeated
application of $F$ will eventually produce the desired elements.

Now for any $\ell$ and any $\bv \in M$, we shall show that
\begin{equation} \label{eq:nabla}
|(\nabla - i)\cdots(\nabla-i-p^\ell+1)\bv| \leq |p^\ell! \bv|.
\end{equation}
(It is this point that is unclear in \cite{bib:tsu3}.)
Of course it suffices to work with $\bv$ equal to one of our chosen elements.
In that case, modulo $(p^\ell)!$, $\nabla$ acts simply as the derivation
$\theta$, so it suffices to show that
\[
|(\theta - i)\cdots(\theta-i-p^\ell+1)\bv| \leq |p^\ell! \bv|.
\]
But this is evident: applying $(\theta-i)\cdots(\theta-i-p^{\ell}+1)$
maps $t^m$ to $(m-i)\cdots(m-i-p^{\ell}+1) t^m$, and the product of
$p^{\ell}$ consecutive integers is always divisible by $p^\ell!$.

With (\ref{eq:nabla}) now proved, we conclude that for $n = j p^\ell +k$
with $0 \leq k < p^{\ell}$,
\begin{align*}
\left| \frac{1}{n!} \nabla^{(n)} \bv\right| &\leq |p^\ell!|^j\, |n!|^{-1} \\
&= |j!k!| \\
&\leq p^{(j + p^\ell - 1)(p-1)} \\
&\leq p^{-c+n/(p^{\ell}(p-1))} < p^{n\epsilon+c},
\end{align*}
as desired.
\end{proof}

\begin{proof}[Proof of Proposition~\ref{prop:changefrob}]
The argument is motivated by the fact that for all $x \in R$,
\[
x^{\sigma_2} = \sum_{n=0}^\infty \frac{1}{n!} (t^{\sigma_2}/t^{\sigma_1}
- 1)^n (\theta (\theta-1)\cdots (\theta-n) x)^{\sigma_1}.
\]
To check this, verify that the right side is a ring endomorphism and notice
that the two sides agree when $x=t$.

Given $M$ an $(F, \nabla)$-crystal over $R$ equipped with $\sigma_1$,
we wish to define a linear map from $M^{\sigma_2}$ to $M^{\sigma_1}$,
with which we can then compose $F$ to get a linear map from $M^{\sigma_2}$
to $M$. We will show that
\begin{equation} \label{eq:changefrob}
\bv \mapsto \sum_{n=0}^\infty \frac{1}{n!} (t^{\sigma_2}/t^{\sigma_1}-1)^n
\nabla(\nabla-1)\cdots(\nabla-n+1) \bv
\end{equation}
is such a map, except that in general it maps $M^{\sigma_2}$ into
$M^{\sigma_1} \otimes_{\calO} \calO[\fp]$. But that will suffice, because
then a suitable isogeny will produce an actual $F$-crystal structure on $M$
equipped with $\sigma_2$. 

First of all, we must show that the series in (\ref{eq:changefrob}) converges
in $\Gamma$. This follows from Lemma~\ref{lem:changefrob2}:
if $|t^{\sigma_2}/t^{\sigma_1}-1| = p^{-a}$, then $a>0$, and applying
Lemma~\ref{lem:changefrob2} with $\epsilon<a$ gives the desired convergence.
This completes the proof in the case $R = \Gamma$.

In case $R = \Gcon$, we must also show that the series in
(\ref{eq:changefrob}) converges to a limit $\bw$ defined over $\Gcon$ and not just
over $\Gamma$. Let $N$ be the matrix through which $\nabla$ acts on our
chosen system of elements, and $c,d>0$ be such that $v_n(N_{ij}) \geq -cn-d$
and $v_n(b) \geq -cn-d$, where $b = t^{\sigma_2}/t^{\sigma_1}-1$.
As in the previous paragraph, we may deduce from the lemma that there exist
constants $e,f >0$ such that all but (at most) the first $en+f$ terms of the series
have absolute value less than $p^{-n}$. Putting
$\nabla^{(i)} = \nabla (\nabla-1) \cdots (\nabla-i+1)$ as in the lemma,
we have
\begin{align*}
v_n(\bw) &\geq \min_{i\leq en+f} \{ v_n(\frac{1}{i!} b^i \nabla^{(i)}\bv ) \} \\
&\geq \min_{i \leq en+f} \{ v_{n+i/(p-1)}(b^i \nabla^{(i)}\bv ) \} \\
&\geq \min_{i \leq en+f} \{ \min_{j \leq n+i/(p-1)} \{v_j(b^i) v_{n+i/(p-1)-j}(\nabla^{(i)} \bv) \} \} \\
&\geq \min_{i \leq en+f} \{ \min_{j \leq n+i/(p-1)} -c(n+i/(p-1)) - (n+i/(p-1))d \} \\
&\geq -(c+d)n - (c+d)(en+f)/(p-1).
\end{align*}
Thus $\bw$ is overconvergent.
\end{proof}

Note that if $M$ is trivial, then the change of Frobenius 
map is none other than the ring endomorphism we gave at the beginning of
the above proof. In other words,
trivial crystals remain trivial under change of Frobenius.

Note that in some cases, notably if $\calO = W$ and $p > 2$, then the
terms in (\ref{eq:changefrob}) are automatically integral. That
means that in those cases, the category of $(F,\nabla)$-crystals up to
\emph{isomorphism} is independent of the definition of Frobenius.
This is also true \emph{a priori} for crystals that come from geometry
(as the crystalline cohomology of varieties), since the structure of
a module on the crystalline site gives rise to compatible $(F,\nabla)$-crystals
for all possible choices of Frobenius, well-defined up to isomorphism.

If a Frobenius $\sigma$ on $R_1$ extends canonically to $R_2$, we
can base-extend $F$-crystals $M$ from $R_1$ to $R_2$, by having $F$
map $(M \otimes_{R_1} R_2)^\sigma
= M^\sigma \otimes_{R_1} R_2^\sigma$ to $M \otimes_{R_1} R_2$ by
the original action of $F$ on the first factor and the inclusion of
$R_2^\sigma$ to $R_2$ on the second factor. Similarly, if $\theta$ also
extends canonically, we can base-extend $(F,\nabla)$-crystals $M$
from $R_1$ to $R_2$.

For $\sigma$ standard, we say
an $F$-crystal over $\Gcon$ is \emph{semistable} if it is isogenous (over $\Gcon$)
to a crystal over $\Omega$. For $\sigma$ arbitrary, we say 
an $(F, \nabla)$-crystal over $\Gcon$ is \emph{semistable}
if there exists a residual uniformizer
$t \in \Gcon$ such that after changing Frobenius to $\sigma_t$,
the crystal becomes isogenous
(over $\Gcon$) to a crystal over $\Omega_t$.
It will follow, once we have established our main result, that
this is equivalent
to the crystal becoming isogenous to a crystal over $\Omega_t$ after
changing Frobenius to $\sigma_t$ for
\emph{every} residual uniformizer $t \in \Gcon$.

We say that
an $(F, \nabla)$-crystal over $\Gcon$
is \emph{potentially semistable} if it becomes stable after making a finite
extension of $k((t))$. More precisely, this extension can always be taken to
be separable. To see this, suppose the
$(F, \nabla)$-crystal $M$ over $\GK_{\con}$ becomes semistable
over $\GL_{\con}$ with $L = K((t^{1/p}))$. Then a standard
Frobenius on $\GL_{\con}$ is also standard on $\GK_{\con}$.
If $\be_1, \dots, \be_n$ span $M$ over $\GL_{\con}$ and $F$ acts
on them through a matrix over $\Omega$, then $\be_1^{\sigma},
\dots, \be_n^{\sigma}$ span $M$ over $\GK_{\con}$ and $F$ also
acts on them through a matrix over $\Omega$. Thus $M$ is already semistable
over $\GK_{\con}$.

\section{Semistability and unipotency}

In this section we prove the following theorem and corollary,
the main results of the paper.
\begin{theorem} \label{mainthm}
Let $M$ be an $F$-crystal over $\Gcon$. Then
$M$ is semistable if and only if $M$ becomes constant over $\Gancon$.
\end{theorem}
\begin{cor}
Let $M$ be an $(F, \nabla)$-crystal over $\Gcon$. Then
$M$ is potentially semistable if and only if $M$ becomes constant as
an $F$-crystal (and unipotent as an $(F, \nabla)$-crystal) over
a finite extension of $\Gancon$.
\end{cor}

\subsection{Dwork's trick}

This section consists of the proof of one half of Theorem~\ref{mainthm},
namely that a semistable $F$-crystal over $\Gcon$ becomes constant
over $\Gancon$. More precisely, we need to show that (with $\sigma = \sigma_t$
standard)
any $F$-crystal over
$\Omega$ becomes constant over $\Oan$. 
This fact is referred to as ``Dwork's trick'' by de~Jong;
we prove a  slightly stronger assertion than his version
\cite[Lemma~6.2]{bib:dej1}, which applies only to $(F, \nabla)$-crystals.

\begin{lemma}[Dwork's trick] \label{lem:dwork}
Every $F$-crystal over $\Oan$ is constant.
\end{lemma}
\begin{proof}
Choose an isobasis $\be_1, \dots, \be_n$ of $M$ such that
$F\be_i \equiv \lambda_i \be_i \pmod{t}$ for some $\lambda_i \in \calO$. Now
define $\be_i^{(n)} = F^n \be_i / \lambda_i^n$ as an element of
$M \otimes_{\calO} \calO[\fp]$. Then
\[
\be_i^{(n+1)} - \be_i^{(n)}
= \lambda_i^{-n} F^n (F\be_i - \lambda_i\be_i) \equiv 0 \pmod{t^{p^n}}.
\]
Thus the sequence $\be_i^{(n)}$ converges $t$-adically in $M \otimes_\Omega
\calO[\fp][[t]]$ to a limit
which we call $\bof_{i}$.

We claim that $\bof_{i}$ is actually defined over $\Oan$.
To show this, let $A$ be the matrix by which $F$ acts on the basis
$\{\be_{1}, \dots, \be_{n}\}$,
and let $\bof_{i} = \be_{i} + \sum_{j} c_{j} \be_{j}$.
Now we have $F\bof_{i}=\lambda_{i}\bof_{i}$, which we
may rewrite as
\begin{equation} \label{dwork1}
\bof_{i} - \be_{i} = \lambda_{i}^{-1} F(\bof_{i} - \be_{i})
+ \lambda_{i}^{-1}F\be_{i}  - \be_{i}.
\end{equation}
Putting $\lambda_{i}^{-1}F\be_{i}  - \be_{i} = \sum_{j} d_{j} \be_{j}$,
we have $d_{j} \in t\Oan$. Now write $A_{jk} = \sum_m A_{jk,m} t^m$,
$c_j = \sum_m c_{j,m} t^m$, and $d_j = \sum_m d_{j,m} t^m$.

Because $d_j \in t\Oan$, for any $c>0$ there exists
$d>0$ such that $|\lambda_i^{-1} A_{jk,n}| \leq p^{-cn-d}$
and $|d_{j,n}| \leq p^{-cn-d}$ for all $n>0$. Moreover, there
exists $N>0$ such that $|\lambda_{i}^{-1} A_{jk},n| \leq p^{-cn}$
for $n \geq N$.
The equation~(\ref{dwork1})
becomes $c_{j} = \sum_{k} A_{jk} c_{k}^{\sigma} + d_{j}$,
giving the estimate
\begin{equation} \label{dwork2}
|c_{j,n}| \geq \max_{k,m} \{
|A_{jk,m} c_{k,(n-m)/p}|, |d_{j,n}| \}.
\end{equation}
A straightforward induction and the fact that $-cn-d \geq -cn-nd$ for
$n >0$ gives the conclusion $|c_{j,n}| \leq p^{-cn-nd}$
for all $n>0$. This bound is not strong enough to give the desired
conclusion, but for large enough $n$ it can be substantially improved.

To be precise, we show that if we put $K = 1+\frac{dp}{c(p-1)}$, then
$|c_{j,n}| \leq p^{-cn-Kd}$ for $n \geq K$. Again, this is
by induction on $n$. For $m<n$, in case $m \leq n-K+1$ we have
\begin{eqnarray*}
-cm -d -c \frac{n-m}{p} - Kd &=& -cm \frac{p-1}{p} - \frac{cn}{p} - (K+1)d \\
&\geq& -c \frac{p-1}{p} \left(n - \frac{dp}{c(p-1)} \right)
- \frac{cn}{p} - (K+1)d \\
&=& -cn - Kd.
\end{eqnarray*}
Otherwise, using the earlier estimate $|c_{j,n}| \leq p^{-cn-dn}$ for $n
< K$, we have
\begin{eqnarray*}
-cm - d - c\frac{n-m}{p} - d\frac{n-m}{p}
&\geq& -cm -d - c(n-m) - d(n-m) \\
&\geq& -cn - Kd.
\end{eqnarray*}
From this and the estimate (\ref{dwork2}), we deduce the desired
inequality by induction.

All that remains is to show that the change of basis matrix from
the $\be_i$ to the $\bof_i$ is invertible. In fact,
the inverse matrix can be constructed as an analogous change of basis
matrix for the dual crystal $M^*$.
\end{proof}

\subsection{Factorization of matrices over $\Gamma$}
The next two sections are devoted to the proof of the remaining half of
Theorem~\ref{mainthm}, that is, that an $F$-crystal over $\Gcon$ which
becomes constant over $\Gancon$ is semistable. Loosely speaking,
we show this by proving that
the change-of-basis matrix over $\Gancon$ can be factored as the product
of a matrix over $\Oan$ times a matrix over $\Gcon$, and changing basis by the
latter gives a presentation of the crystal over $\Omega$.

The rank 1 case of the following assertion can be found in Zannier \cite{bib:zan}.
\begin{lemma} \label{lem:factor}
Let $U = \sum_n U_n t^n$ be a matrix over $\Gamma$ such that $U \equiv 1
\pmod{\pi}$. Then there exists a unique pair $(P, N)$ of matrices over
$\Gamma$ of the form $P = 1 + \sum_{n=1}^\infty A_n t^n$ and $N = C +
\sum_{n=1}^\infty B_n t^{-n}$, with $C$ a matrix over $\calO$,
such that $U = PN$.
Moreover, if $|U_n| < \epsilon$ for $n > a$ (resp. $n < -b$),
then $|A_n| < \epsilon$ for $n > a$ (resp. $|B_n| < \epsilon$ for $n < -b$).
\end{lemma}
\begin{proof}
We define convergent sequences $\{P_{n}\}, \{N_{n}\}$ such that $P_{n}
N_{n} \equiv U \pmod{\pi^{n}}$. We start with $P_{1}=N_{1}=1$. To
define $P_{n+1}$ and $N_{n+1}$, write
\[
P_{n}^{-1} U N_{n}^{-1} = \pi^{n}(A + B),
\]
with $A = \sum_{i=1}^{\infty} A_{i} t^{i}$ and $B =
\sum_{i=0}^{\infty} B_{i} t^{-i}$, and put $P_{n+1} = P_n (1+\pi^{n} A)$
and $N_{n+1} = (1+\pi^{n}B) N_n$. Then the sequences $\{P_n\}$ and 
$\{N_n\}$ both converge $p$-adically to the desired $P$ and $N$.

To establish the final assertion, write $N^{-1} = C^{-1} + 
\sum_{n=1}^\infty D_n t^{-n}$, and note that since $P = UN^{-1}$,
\[
A_n = U_n C^{-1} + \sum_{m=1}^{\infty} U_{n+m} D_{m}
\]
and so if $n > a$, then $|A_n| < \epsilon$; the proof that $|B_n| < 
\epsilon$
for $n < -b$ is similar.
\end{proof}

Recall that an \emph{elementary matrix} is a matrix obtained from the
identity by adding a multiple of one row to another, swapping two rows, or
multiplying one row by a unit. The following lemma can be deduced easily from
\cite[Theorem~III.7.9]{bib:lang}.
\begin{lemma} \label{lem:rowop}
Let $R$ be a principal ideal domain.
\begin{enumerate}
\item[(a)]
Let $A$ be an invertible $n \times n$ matrix over $R$.
Then $A$ is the product of elementary matrices.
\item[(b)]
    Let $A$ be an $n \times n$ matrix over $R$
with
    determinant 0. Then there exists an invertible matrix $B$
such that $AB$ has zeroes in its first column.
\end{enumerate}
\end{lemma}

\begin{lemma} \label{lem:rowop2}
Let $V$ be a matrix over $\Gamma$ with nonzero determinant,
such that for each $n \in \NN$, $V$ is congruent modulo $\pi^n$ to
a matrix over $\calO[t,t^{-1}]$. Then
there exists an invertible matrix $B$ such that:
\begin{enumerate}
    \item[(a)] $B$ and $B^{-1}$ have entries in $\calO[t, t^{-1}, \fp]$;
    \item[(b)] $VB$ has integral entries and $|VB-1| < 1$.
\end{enumerate}
\end{lemma}
\begin{proof}
    We first prove the statement with a weaker form of (b), namely that
$VB$ has integral entries and determinant not divisible by $\pi$.
For this, we may induct on the valuation of $\det(V)$. The case where this
    valuation is zero is trivial, so we assume $\det(V) \equiv 0 \pmod{\pi}$.
    By Lemma~\ref{lem:rowop}, there exists an invertible matrix $B$
    over $\Gamma$ satisfying (a), and such that the entries of the first
    column of $VB$ are divisible by $\pi$. (Write $B$ in the
    formulation of the lemma as a product of elementary matrices over $k[t, t^{-1}]$ and
    lift each to an elementary matrix over $\calO[t,t^{-1}]$.
    Then $B^{-1}$ is also a product of
    such matrices.) We now can multiply $B$ by a diagonal matrix on
        the right so as to divide the entries in the first column of $VB$ by $\pi$. This
    reduces the valuation of $\det(V)$ while maintaining the integrality
    of the entries, so application of the induction hypothesis completes
    the proof of the weaker assertion.

To prove the original assertion, it suffices to note that by
Lemma~\ref{lem:rowop} again, the reduction of $VB$ is the product of
elementary matrices, so again it can be lifted to a product $C$ of
elementary matrices so that $C$ and $C^{-1}$ have entries which are
finite sums of powers of $t$.  Replacing $B$ with $BC$ gives the
desired result.
\end{proof}

\subsection{Factorization of matrices over $\Gancon$}

To prove Theorem~\ref{mainthm}, we need a lemma to the effect that
given a matrix over $\Gancon$, one can ``factor off the part not
defined over $\Gamma$''. We will deduce such a lemma using the results
of the previous section, by ``tilting'' a matrix over $\Gancon$
to put its entries into $\Gcon$.

In passing, we note that Berger \cite{berger} has used the following lemma
to give a simplified proof of Colmez's theorem that absolutely crystalline
representations are of finite height (conjectured by Fontaine).

\begin{lemma} \label{lem:anfact1}
Let $U$ be a matrix over $\Gancon$ with nonzero determinant. Then there exist
a matrix $V = 1 + \sum_{n=1}^\infty V_n t^n$ over $\Oan$
and a matrix $W$ over $\Gcon[\fp]$ such that $U = VW$.
\end{lemma}
\begin{proof}
Put $U = \sum_{n=-\infty}^\infty U_n t^n$.
Let $c$ be a positive rational number strictly less than
$\liminf_{n \to +\infty} \{(\log_p |U_n|)/n\}$. For the moment, we enlarge the
ring of scalars by replacing $\calO$ with $\calO[p^c]$.

The operation of tilting consists of replacing
a matrix $X = \sum_i X_i t^i$ with the new matrix $\tilde{X} = 
\sum_i X_i p^{ci} t^i$. The tilted matrix $\tilde{U}$
has coefficients in $\Gcon[\fp]$ and, for each $n$,
is congruent modulo $p^n$ to a \emph{finite} sum
$\sum_{i} W_{i} t^{i}$. There is no loss of generality in assuming that 
$\tilde{U}$ has integral entries (by multiplying $U$ by an appropriate scalar
matrix, which in the end can be divided from $U$ and $W$).

Apply Lemma~\ref{lem:rowop2} with $V = \tilde{U}$, and
let $\tilde{B}$ be the matrix $B$ in
the conclusion of the lemma.
Then Lemma~\ref{lem:factor} gives a decomposition 
$\tilde{U}\tilde{B} = \tilde{P} \tilde{N}$, where $\tilde{P}$
has only positive powers of $t$ (and constant term 1),
and $\tilde{N}$
has only negative powers of $t$ (and invertible constant term).

We now wish to untilt $\tilde{B}, \tilde{P}, \tilde{N}$ and
conclude that we still have the decomposition $UB = PN$. For this, it suffices
to show that $B, P, N$ have entries in $\Gancon$.
This is obvious
for $B$ because it has entries in $\calO[t,t^{-1},\fp]$,
and for $N$ because untilting a matrix with only negative
powers of $t$ only decreases the absolute values of its coefficients.

As for $P$, note that for each $\epsilon > 0$,
there exists $d$ such that $|(UB)_n| < p^{\epsilon n+d}$ for $n$ large enough,
and so $|(\tilde U \tilde B)_n|
< p^{-cn+\epsilon n + d}$ for $n$ large enough.
Applying the final assertion of Lemma~\ref{lem:factor}, we get that
$|\tilde{P}_n| < p^{-cn+\epsilon n + d}$ for $n$ large enough, so
$|P_n| < p^{\epsilon n + d}$ for $n$ large enough. We conclude that $P$
has entries in $\Oan$.

Set $V = P$ and $W = NB^{-1}$. We now have the desired
factorization $U=VW$, except that we have enlarged $\calO$ and the assertion
of the lemma does not permit such an enlargement.
On the other hand, if $m_1, \dots, m_i$ is a basis for $\calO[p^c]$
over $\calO$, with $m_1 = 1$, and $V = \sum_{j=1}^i m_j V_j$, then
$U^{-1}V = \sum_{j=1}^i m_j (U^{-1}V_j)$. Since $U^{-1}V$ has entries in
$\Gcon[\fp] \otimes_{\calO}
\calO[p^c]$, $U^{-1}V_j$ has entries in $\Gcon[\fp]$ for $j=1, \dots, i$.
In particular, $W_1^{-1} = U^{-1}V_1$ has entries in $\Gcon[\fp]$, and
$U = V_1W_1$ is a decomposition of the desired form.
\end{proof}

We wish to refine the decomposition given by the previous lemma under
the additional assumption that $U$ is invertible over $\Gancon$. To do so,
we first identify the units in $\Gancon$.
\begin{lemma} \label{lem:units}
For $x = \sum_i x_i t^i$ a nonzero element of $\Gancon$, $x$ is invertible
if and only if $|x_i|$ is bounded above.
\end{lemma}
\begin{proof}
Clearly if $|x_i|$ is bounded above, then $x \in \Gcon[\fp]$ is invertible.
Conversely, suppose $xy = 1$ but $|x_i|$ is not bounded above.
Choose $c>0$ such that $|x_i|p^{ci}$ and $|y_i| p^{ci}$ are bounded
above for $i \in \QQ$. Let $i_1, i_2, \dots$ be the sequence of indices $i$ such that
$|x_j| < |x_i|$ for $j<i$, which by assumption is infinite.
Put
\[
s_k = \frac{v_p(x_{i_{k}}) - v_p(x_{i_{k+1}})}{i_{k+1} - i_k};
\]
then $s_k \to 0$ as $k \to \infty$. So we may choose $k$ such that $s_k < c$.
Now tilt $x$ and $y$ to obtain $\tilde{x} = \sum_i x_i p^{ci} t^i$
and $\tilde{y} = \sum_i y_i p^{ci} t^i$. 
On one hand, we have $\tilde{x}\tilde{y}
= 1$. On the other hand, 
we may assume without loss of generality that
$\tilde{x}$ and $\tilde{y}$ are elements of $\Gcon$ with nonzero reduction;
then $\tilde{x}$ and $\tilde{y}$ are congruent modulo $\pi$ to polynomials
which are not both constant. (By construction, $\tilde{x}_{i_k}$ and
$\tilde{x}_{i_{k+1}}$ are nonzero modulo $\pi$.) Thus their product
is congruent to a nonconstant polynomial modulo $\pi$, contradiction.
\end{proof}
\begin{cor}
Suppose $x \in \Oan$ becomes invertible in $\Gancon$. Then there exists
a polynomial $y \in \calO[t,\fp]$ such that $x/y$ is invertible in $\Oan$.
\end{cor}
\begin{proof}
By Lemma~\ref{lem:units}, we have $x \in \Oan \cap \Gcon[\fp] = \Omega[\fp]$.
The desired result then follows immediately from the Weierstrass preparation theorem
\cite[Theorem~IV.9.2]{bib:lang}.
\end{proof}

\begin{lemma} \label{lem:anfact2}
Let $A = 1 + \sum_{n=1}^\infty A_n t^n$ be an $m \times m$ matrix over
$\Oan$ which becomes invertible over $\Gancon$. Then there exist a matrix
$B = 1 + \sum_{n=1}^\infty B_n t^n$ which is invertible over $\Oan$
and a matrix $C$ over $\Gcon[\fp]$ such that $A = BC$.
\end{lemma}
\begin{proof}
The proof will resemble that of Lemma~\ref{lem:rowop2}.
By the previous corollary, $\det A$ has finitely many zeroes in the formal
unit disc, and we induct on the number of these zeroes, counted with
multiplicities.

If $\det A$ has no zeroes in the disc, its
inverse is rigid analytic (again by the previous corollary)
and so we may use the trivial factorization
$A = A \cdot 1$. Otherwise, let $r$ be a zero of $\det A$,
let $P(t)$ be the minimal polynomial of $r$ over $\calO$,
and let $D$ be the reduction of $A$ modulo $P(t)$.
By assumption, $\det D = 0$, so
the rank of $D$ is less than $m$. Therefore there exists an invertible
matrix
$E$ over $\calO[t,\fp]$ such that $DE$, and likewise $AE$,
has its first column identically
zero modulo $P(t)$.
Now let $F$ be the diagonal matrix with $P(t)$ as its first entry
and 1 in its other diagonal positions. The matrix $AEF^{-1}$ now satisfies
the same hypotheses as $A$, but $\det AEF^{-1}$ has fewer zeroes than
does $\det A$. By the induction hypothesis, we have $AE F^{-1} = BC$,
with $B$ invertible over $\Oan$ and having constant term 1, and $C$
defined over $\Gcon[\fp]$. We now factor $A$ as $B(CF^{-1}E)$,
and the two terms again have the desired properties.
\end{proof}

\begin{cor} \label{cor:anfact}
Let $U$ be an invertible matrix over $\Gancon$. Then there exist
a matrix $V = 1 + \sum_{n=1}^\infty V_n t^n$ which is invertible
over $\Oan$
and a matrix $W$ over $\Gcon[\fp]$ such that $U = VW$.
\end{cor}
\begin{proof}
By Lemma~\ref{lem:anfact1}, we have a decomposition $U = VW$, with $V$
having entries in $\Oan$ (and constant term 1) and $W$ having entries
in $\Gcon[\fp]$. In particular, $V = UW^{-1}$ is invertible in $\Gancon$,
so Lemma~\ref{lem:anfact2} allows us to write $V$ as $XY$, with $X$
invertible over $\Oan$ (and having constant term 1), and $Y$ having
entries in $\Gcon[\fp]$. The decomposition $U = X(YW)$ now has
the desired form.
\end{proof}

We can now complete the proof of Theorem~\ref{mainthm}.
Suppose $M$ is an $F$-crystal over $\Gcon$ which
becomes constant over $\Gancon$. That means that if $A$
is the matrix representing the action on Frobenius on an isobasis of $M$
over $\Gcon$, then there exists an invertible matrix $U$ over
$\Gancon$ and a matrix $D$ over $\calO$ such that $A =
U^{-1} D U^\sigma$. By Corollary~\ref{cor:anfact}, we may write $U = VW$,
with $V$ invertible over $\Oan$ and $W$ defined over $\Gcon[\fp]$.
Changing basis by $W$ gives a new isobasis on which
Frobenius acts
via the invertible matrix $V^{-1}DV^\sigma$ over $\Oan \cap \Gcon[\fp] = \Omega[\fp]$.
Thus $M$ is isogenous to a crystal defined over $\Omega$,
that is, $M$ is semistable, as desired.

As noted earlier, the equivalence of semistability and unipotency implies
that semistability, which is defined by imposing a condition on the change
to a single standard Frobenius, can be checked by changing to any standard
Frobenius. It also implies that the conjecture that every overconvergent
$(F, \nabla)$-crystal is quasi-unipotent \cite[Section~10.1]{bib:crew2},
is equivalent to the following
conjecture.
\begin{conj}
Every $(F, \nabla)$-crystal over $\Gcon$ is potentially semistable.
\end{conj}

This conjecture is known for isoclinic crystals by a theorem of Tsuzuki
\cite{bib:tsu1}. Additionally, the conjecture holds for crystals ``of
geometric origin'' (in a sense to be made precise in a subsequent paper).
In some cases, it can even be established by explicit computation
together with Theorem~\ref{mainthm}.
For example, Tsuzuki \cite{bib:tsu3} showed
that the Bessel crystal (constructed
originally by Dwork) is quasi-unipotent by a direct computation; consequently,
we may conclude
that it is potentially semistable without explicitly computing an
appropriate change of basis.

\section*{Acknowledgments}
This work is based on the author's doctoral dissertation \cite{bib:methesis},
written under the supervision of Johan de~Jong. This paper was written while
the author was supported by a Clay Mathematics Institute
Liftoffs fellowship and by a National Science Foundation Postdoctoral
Fellowship. The author also thanks Laurent Berger for pointing out
his preprint \cite{berger}.


\end{document}